\def \CC {\mathbb C}

\def \epsilon{\varepsilon}

\def \S  {{\mathcal S}}

\def \ga {\gamma}

\renewcommand{\S}{{\mathcal S}}

\documentclass[12pt,reqno]{amsart}

\usepackage{amsmath,amssymb}

\usepackage{url} 

\numberwithin{equation}{section}

\begin{document}


\title[]{A note on the invariants of the $L$-functions} 

\author[]{J.KACZOROWSKI \lowercase{and} A.PERELLI}
\maketitle

{\bf Abstract.} We explain the exact meaning of a statement we made in a previous paper on invariants, namely that a complex-valued function of the data of the functional equation of an $L$-function is an invariant if and only if it is stable under the multiplication and factorial formulae. To this end, we show that every invariant has a so-called rational extension, having some desired invariance properties. The existence of such an extension, which enables to express formally the above heuristic concept, is not apparent and its construction is the main novelty of the paper.

\smallskip
{\bf Mathematics Subject Classification (2010):} 11M41

\smallskip
{\bf Keywords:} Selberg class; gamma factors; invariants.

\vskip.5cm
\section{Introduction}

\smallskip
Every element $F$ of the extended Selberg class, denoted by $\S^{\sharp}$, satisfies a Riemann-type functional equation of the form
\[ 
\gamma(s) F(s) = \omega \overline{\gamma}(1-s)\overline{F}(1-s),
\]
where $\omega\in {\mathbb C}$, $|\omega|=1$ and
\begin{equation}
\label{eq:gamma}
 \gamma(s) = Q^s\prod_{j=1}^r\Gamma(\lambda_j s+\mu_j)
 \end{equation}
with  $Q>0$, $\lambda_j>0$ and $\Re(\mu_j)\geq 0$  for $j=1,\ldots,r$. Here $\Gamma(s)$ denotes the familiar Euler gamma function and $\overline{f}(s)=\overline{f(\overline{s})}$. The numbers and vectors
\begin{equation}
\label{eq:inv}
 \omega, Q, \boldsymbol{\lambda}=(\lambda_1,\ldots,\lambda_r), \boldsymbol{ \mu}=(\mu_1,\ldots,\mu_r)
 \end{equation}
are called the data of the functional equation. It is well-known that the $\gamma$-factor $\gamma(s)$, as a function of a complex variable $s$, is uniquely determined by $F$ up to a non-zero multiplicative constant; see Conrey-Ghosh \cite{Co-Gh/1993}. Nevertheless, the functional equation of a given $F\in\S^\sharp$ may have different forms with different data sets. This is due to the functional equations satisfied by the Euler function, in particular, the factorial formula
\[ 
s\Gamma(s)=\Gamma(s+1)
\]
and the Legendre-Gauss multiplication formula
\[
\Gamma(s)=m^{s-1/2}(2\pi)^{(1-m)/2}\prod_{k=0}^{m-1}\Gamma\left(\frac{s+k}{m}\right).
\]
This leads to the notion of invariant, namely a function
\[ 
I(\omega, Q,  \boldsymbol{\lambda}, \boldsymbol{ \mu})
\]
depending on the data \eqref{eq:inv} whose value depends only on $F$, independently of the particular form of the $\gamma$-factor.

\smallskip
In \cite{Ka-Pe/2000}, we proved that any $\gamma$-factor of $F\in\S^\sharp$ can be transformed to a complex multiple of any other $\gamma$-factor of $F$ by repeated applications of the above two formulae. As a consequence, we sometimes used the following heuristic assertion:
\begin{equation}
\label{prop}
\begin{split}
&\text{{\it $I(\omega, Q,  \boldsymbol{\lambda}, \boldsymbol{ \mu})$ is an invariant if and only if it is}} \\
& \hskip-.5cm\text{{\it stable under the multiplication and factorial formulae,}}
\end{split}
\end{equation}
where ``stable'' meant, generically, ``suitably invariant''. The authors were then frequently asked about the precise meaning of this assertion. While the meaning of stability under the multiplication formula is obvious and was rigorously defined in \cite{Ka-Pe/2000}, the problem is with the application of the factorial formula. Indeed, applying such a formula to $\Gamma(\lambda_1s+\mu_1)$ in \eqref{eq:gamma} changes $\gamma(s)$ to
\[ 
(\lambda_1s+\mu_1-1) \Gamma(\lambda_1 s+\mu_1-1) Q^s \prod_{j=2}^r\Gamma(\lambda_j s+\mu_j),
\]
which is no longer of the form \eqref{eq:gamma}. In general, repeated applications of the factorial formula  transform a $\gamma$-factor into 
\begin{equation}
\label{eq:Rgamma}
 R(s)\gamma(s),
 \end{equation}
where $R(s)$ is a rational function and $\gamma(s)$ is of the form \eqref{eq:gamma}. Formally, for such expressions, $I(\omega, Q,  \boldsymbol{\lambda},\boldsymbol{ \mu})$ is not defined; hence the meaning of stability is not clear. Of course the problem disappears if the factorial formula is not used in transforming a $\gamma$-factor to a multiple of another $\gamma$-factor. According to Theorem 1 in \cite{Ka-Pe/2000}, this holds when the $\gamma$-class number $h_F$ of $F$ is $1$ (meaning that in \eqref{eq:gamma} we have  $\lambda_i\slash\lambda_j\in {\mathbb Q}$ for $i,j=1,\ldots,r$) or if $F$ is reduced (meaning that $0\leq \Re(\mu_j)<1$ for $j=1,\ldots,r$). Another instance where an application of the factorial formula does not create the above-mentioned problems occurs when $I(\omega, Q,  \boldsymbol{\lambda}, \boldsymbol{ \mu})$ depends only on $\omega,Q,\boldsymbol{\lambda}$ and $\Im(\boldsymbol{\mu}):= (\Im(\mu_1),\ldots,\Im(\mu_r))$. Indeed, applications of the factorial formula may change the real parts of the $\mu_j$'s only, so in such cases the appearance of the rational factor in \eqref{eq:Rgamma} can be ignored. For instance, the factorial formula can be avoided, or the presence of the  rational factor can be ignored, in the case of the following well-known invariants: the degree
\[ 
d_F= 2\sum_{j=1}^r \lambda_j,
\]
the conductor
\[ 
q_F = (2\pi)^{d_F}Q^2\prod_{j=1}^r \lambda_j^{2\lambda_j}
\]
and the root number
\[ 
\omega_F= \omega\prod_{j=1}^r \lambda_j^{-2i\Im(\mu_j)}.
\]

\smallskip
In contrast, there are problems in the case of the important $H$-invariants, defined in \cite{Ka-Pe/2002} for every integer $n\geq0$ as
\begin{equation}
\label{eq:H}
 H_F(n)=2\sum_{j=1}^r \frac{B_n(\mu_j)}{\lambda_j^{n-1}},
 \end{equation}
where $B_n(x)$ denotes the $n$-th Bernoulli polynomial. Checking the stability of $H_F(n)$  under the multiplication formula is not difficult, see e.g. Section 3.3 of \cite{Kac/2006}. We only remark that it is performed using the following well-known property of the Bernoulli polynomials
\[
B_n(mx)=m^{n-1}\sum_{j=0}^{m-1} B_n(x+\frac{j}{m}) \hskip1.5cm (n\geq 0,\ m\geq 1).
\]
However, the treatment of the stability under the factorial formula in \cite{Kac/2006} is not entirely satisfactory, as it treats a special case where the formula is applied to a pair of factors $\Gamma(\lambda s+\mu)$ and $\Gamma(\lambda's+\mu')$  satisfying the extra consistency condition
\[
\frac{\mu-1}{\lambda}=\frac{\mu'}{\lambda'}.
\]
It is not evident that the algorithm transforming a $\gamma$-factor to another can be arranged so that applications of the factorial formula are always performed in the above way. Hence this needs further explanation, otherwise the argument cannot be regarded as general. We stress that these comments apply only to the method of proof used in \cite{Kac/2006} to check the stability of $H_F(n)$ under multiplication and factorial formulae. Indeed, there are other proofs of the fact that the $H$-invariants are actually invariants, based on different principles; see \cite{Ka-Pe/2002}.

\smallskip
To clarify the situation with $H_F(n)$, and with the general case as well, we introduce the following definitions. Let $R(s)$ denote a generic not identically vanishing rational function, $\gamma(s)$ be a $\gamma$-factor as in \eqref{eq:gamma} and $I(\omega, Q,  \boldsymbol{\lambda}, \boldsymbol{ \mu})$ be a function of the data of $\gamma(s)$. Let now $I^*(R,\omega, Q,  \boldsymbol{\lambda}, \boldsymbol{ \mu})$ be a complex-valued function depending on the parameters involved in $R(s)$ and on the data of $\gamma(s)$. In other words, we may think of $I^*$ as a function of the data of a product of type \eqref{eq:Rgamma}. We say that  $I^*$ is a {\it rational extension} of  $I$ if for every {\it constant} rational function $R$ we have 
\[
I^*(R,\omega, Q,  \boldsymbol{\lambda}, \boldsymbol{ \mu}) = I(\omega, Q,  \boldsymbol{\lambda}, \boldsymbol{ \mu}).
\]
An application of the multiplication formula transforms $R(s)\gamma(s)$ to $R(s)c\gamma'(s)$, with certain $c\in\CC\setminus\{0\}$ and $\gamma'(s)$ of type \eqref{eq:gamma}, while the factorial formula transforms it to $R'(s)\gamma'(s)$, with a rational function $R'(s)\not\equiv 0$ and $\gamma'(s)$ again of type \eqref{eq:gamma}. Hence in both cases $R(s)\gamma(s)$ is transformed to a product of a similar form. Accordingly, we say that $I^*$ is {\it stable} under the multiplication and factorial formulae if, with obvious notation,
\begin{equation}
\label{mult}
I^*(R,\omega, Q,  \boldsymbol{\lambda}, \boldsymbol{ \mu}) = I^*(R,\omega\overline{c}/c, Q',  \boldsymbol{\lambda'}, \boldsymbol{ \mu'})
\end{equation}
and 
\begin{equation}
\label{fact}
I^*(R,\omega, Q,  \boldsymbol{\lambda}, \boldsymbol{ \mu}) = I^*(R',\omega, Q',  \boldsymbol{\lambda'}, \boldsymbol{ \mu'}),
\end{equation}
respectively. We can finally state the precise version of the heuristic statement \eqref{prop}.

\medskip
{\bf Theorem. }{\it A function $I(\omega, Q,  \boldsymbol{\lambda}, \boldsymbol{ \mu})$ of the data is an invariant if and only if it has a rational extension which is stable under the multiplication and factorial formulae.}

\medskip
The proof contains an explicit construction of the rational extension, and this is the main novelty of the paper.

 \medskip
{\bf Acknowledgements.}  This research was partially supported by the Istituto Nazionale di Alta Matematica, by the MIUR grant PRIN-2017 {\sl ``Geometric, algebraic and analytic methods in arithmetic''} and by grant 2021/41/B/ST1/00241  from the National Science Centre, Poland.

\section{Rational extension of $H$-invariants}

\smallskip
For a product of type \eqref{eq:Rgamma}, say
\begin{equation}
\label{eq:Rgamma1}
R(s)\gamma(s)=\kappa\frac{\prod_{j=1}^l(s-\alpha_j)}{\prod_{j=1}^m(s-\beta_j)} Q^s \prod_{j=1}^r\Gamma(\lambda_js+\mu_j),
\end{equation}
for $n\geq 1$ we write synthetically
\begin{equation}
\label{eq:Hstar}
H^*(n)= H^*(n; R\gamma) := 2\sum_{j=1}^r\frac{B_n(\mu_j)}{\lambda_j^{n-1}} + (-1)^n 2n\left(\sum_{j=1}^m\beta_j^{n-1} - \sum_{j=1}^l \alpha_j^{n-1}\right)
\end{equation}
with the usual convention that empty products equal $1$ and empty sums equal $0$. It is obvious that $H^*(n)$ is a rational extension of the $H$-invariants $H_F(n)$ defined in \eqref{eq:H}.

\smallskip
Let us check now that $H^*(n)$ is stable under the multiplication and factorial formulae. Since application of the former does not change $R(s)$, in this case the stability of $H^*(n)$ is checked exactly as for $H_F(n)$, see Section 3.3 of \cite{Kac/2006}. The factorial formula can be applied to a single $\Gamma$-factor in \eqref{eq:gamma} in two ways, which in \cite{Ka-Pe/2000} we called ``expanding'' and ``contracting''. Without losing generality we focus on $\Gamma(\lambda_1s+\mu_1)$. The expanding process changes it to $(\lambda_1s+\mu_1-1)\Gamma(\lambda_1s+\mu_1-1)$, whereas the contracting process to $(\lambda_1s+\mu_1)^{-1}\Gamma(\lambda_1s+\mu_1+1)$. Accordingly, $R(s)\gamma(s)$ in \eqref{eq:Rgamma1} is transformed to $R'(s)\gamma'(s)$, where
\[
R'(s) = \begin{cases}
(\lambda_1s+\mu_1-1)R(s)& \text{in the expanding case}\\
(\lambda_1s+\mu_1)^{-1}R(s) &\text{in the contracting case}
\end{cases}
\]
and
\[
\gamma'(s) = \begin{cases}
\Gamma(\lambda_1s+\mu_1-1)\prod_{j=2}^r\Gamma(\lambda_js+\mu_j)& \text{in the expanding case}\\
\Gamma(\lambda_1s+\mu_1+1)\prod_{j=2}^r\Gamma(\lambda_js+\mu_j) &\text{in the contracting case}.
\end{cases}
\] 
Hence in the expanding case we have
\[
\begin{split}
H^*(n,R'\gamma') &= 2\frac{B_n(\mu_1-1)}{\lambda_1^{n-1}} + 2\sum_{j=2}^r \frac{B_n(\mu_j)}{\lambda_j^{n-1}} \\ 
&\hskip1.5cm + (-1)^n 2n\left(\sum_{j=1}^m\beta_j^{n-1} - \frac{(1-\mu_1)^{n-1}}{\lambda_1^{n-1}} - \sum_{j=1}^l \alpha_j^{n-1}\right) \\
&= H^*(n,R\gamma) + \frac{2}{\lambda_1^{n-1}} (B_n(\mu_1-1) 
- B_n(\mu_1)+ n(\mu_1-1)^{n-1}) \\
&= H^*(n,R\gamma),
\end{split}
\]
where in the last step we applied the well known functional equation
\[
B_n(x+1)=B_n(x) + n x^{n-1}.
\]
Since the contracting case can be treated exactly in the same way, the stability of $H^*(n)$ under the factorial formula follows.

\section{Proof of the Theorem}

\smallskip
Suppose that $\ga(s),\ga'(s)$ are two $\ga$-factors of $F$ and $I(\omega, Q,  \boldsymbol{\lambda}, \boldsymbol{ \mu})$ is a function of the data with rational extension $I^*(R,\omega, Q,  \boldsymbol{\lambda}, \boldsymbol{ \mu})$ stable under the multiplication and factorial formulae. Then $\ga(s)$ can be transformed by means of repeated applications of such formulae to $c\ga'(s)$, and the $\omega$-datum of $\ga'(s)$ is therefore $\omega'=\omega\overline{c}/c$. Moreover we have 
\[
I(\omega, Q,  \boldsymbol{\lambda}, \boldsymbol{ \mu}) = I^*(1,\omega, Q,  \boldsymbol{\lambda}, \boldsymbol{ \mu})
\] 
and, after repeated application of the two formulae,  thanks to \eqref{mult} and \eqref{fact} also
\[
I^*(1,\omega, Q,  \boldsymbol{\lambda}, \boldsymbol{ \mu}) = I^*(c,\omega', Q',  \boldsymbol{\lambda'}, \boldsymbol{ \mu'}).
\]
But $ I^*(c,\omega', Q',  \boldsymbol{\lambda'}, \boldsymbol{ \mu'}) =  I(\omega', Q',  \boldsymbol{\lambda'}, \boldsymbol{ \mu'})$, thus $I(\omega, Q,  \boldsymbol{\lambda}, \boldsymbol{ \mu}) = I(\omega', Q',  \boldsymbol{\lambda'}, \boldsymbol{ \mu'})$ and hence $I(\omega, Q,  \boldsymbol{\lambda}, \boldsymbol{ \mu}) $ is an invariant.

\smallskip
Viceversa, suppose that  $I(\omega, Q,  \boldsymbol{\lambda}, \boldsymbol{ \mu})$ is an invariant. Thus its values are independent of the particular choice of $\ga(s)$ and $\omega$ in the functional equation of $F\in\S^\sharp$. On the other hand, we know that the invariants $\omega_F, q_F$ and $H_F(n)$, $n\geq 0$, determine the functional equation of $F$; see Theorem 1 of \cite{Ka-Pe/2002}. Thus 
\[ 
I(\omega, Q,  \boldsymbol{\lambda}, \boldsymbol{ \mu})  = \tilde{I(}\omega_F, q_F,  (H_F(n))_{n=0}^{\infty})
\]
for a certain function $\tilde{I}$. Given $(\omega,Q,  \boldsymbol{\lambda}, \boldsymbol{ \mu})$ and a rational function $R(s)$ we write
\[ 
\gamma_{Q,  \boldsymbol{\lambda}, \boldsymbol{ \mu}}(s)= Q^s\prod_{j=1}^r\Gamma(\lambda_j s+\mu_j)
\]
and
\[ 
I^*(R, \omega, Q,  \boldsymbol{\lambda}, \boldsymbol{ \mu})  = \tilde{I}\left(\omega\prod_{j=1}^r \lambda_j^{-2i\Im(\mu_j)}, 
(2\pi)^{d_{ \boldsymbol{\lambda}}}Q^2\prod_{j=1}^r \lambda_j^{2\lambda_j},  
(H^*(n, R\gamma_{Q,  \boldsymbol{\lambda}, \boldsymbol{ \mu}})_{n=0}^{\infty}\right),
\]
 where $d_{ \boldsymbol{\lambda}}= 2\sum_{j=1}^r\lambda_j$ and $H^*(n, R\gamma_{Q,  \boldsymbol{\lambda}, \boldsymbol{ \mu}})$ is defined by \eqref{eq:Hstar} for $n\geq1$, while for $n=0$ we put $H^*(0,R\gamma_{Q,  \boldsymbol{\lambda}, \boldsymbol{ \mu}}):=d_{ \boldsymbol{\lambda}}$. This is obviously a rational extension of $I$, and is stable under the multiplication and factorial formulae thanks to the known invariance properties of $\omega_F$, $q_F$ and $H^*(n)$. The proof is now complete, and an algorithm for the construction of a rational extension is also provided. \qed

\bigskip
\bigskip

\ifx\undefined\bysame{poly}.
\newcommand{\bysame}{\leavevmode\hbox to3em{\hrulefill}\ ,}
\fi

\bigskip
\bigskip
\noindent
Jerzy Kaczorowski, Faculty of Mathematics and Computer Science, A.Mickiewicz University, 61-614 Pozna\'n, Poland and Institute of Mathematics of the Polish Academy of Sciences, 00-956 Warsaw, Poland. e-mail: \url{kjerzy@amu.edu.pl}

\medskip
\noindent
Alberto Perelli, Dipartimento di Matematica, Universit\`a di Genova, via Dodecaneso 35, 16146 Genova, Italy. e-mail: \url{perelli@dima.unige.it}

\end{document}